\documentclass[a4paper,fleqn]{cas-dc}

\usepackage[numbers]{natbib}
\usepackage{float}
\usepackage{booktabs}
\usepackage{algorithm}
\usepackage{graphicx}
\usepackage[caption=false]{subfig}
\usepackage{caption}
\usepackage{lipsum}
\usepackage{algpseudocode}
\usepackage{amsmath}
\usepackage{amssymb}
\usepackage{booktabs}
\usepackage{arydshln}
\usepackage{lipsum}
\usepackage{natbib}
\usepackage{comment}
\usepackage{color}
\usepackage{cases}
\usepackage{seqsplit}
\usepackage{mwe}
\usepackage{multirow}
\usepackage{epstopdf}
\usepackage{csquotes}
\usepackage{amsthm}
\usepackage{adjustbox}
\usepackage{colortbl}
\usepackage{lineno,hyperref}
\usepackage{babel}
\usepackage{multicol}
\usepackage{multirow}
\usepackage{siunitx}

\def\tsc#1{\csdef{#1}{\textsc{\lowercase{#1}}\xspace}}
\tsc{WGM}
\tsc{QE}
\tsc{EP}
\tsc{PMS}
\tsc{BEC}
\tsc{DE}
\begin{document}  \sloppy
% \let\WriteBookmarks\relax
% \def\floatpagepagefraction{1}
% \def\textpagefraction{.001}
% \shorttitle{}
% \shortauthors{Meysam Mahjoob}
% \author{Meysam Mahjoob}
% \author{Seyed Sajjad Fazeli}
% \author{Soodabeh Milanloui}
% \author{Leyla Saadat Tavasoli}
% \title [mode = title]{A Modified Adaptive Genetic Algorithm for Multi-product Multi-period Inventory Routing Problem}       
% \author{Meysam Mahjoob}

% \begin{abstract}
% Recent developments in urbanization and e-commerce have pushed businesses to deploy efficient systems to decrease their supply chain cost. Vendor Managed Inventory (VMI) is one of the most widely used strategies to effectively manage supply chains with multiple parties. VMI implementation asks for solving the Inventory Routing Problem (IRP). This study considers a multi-product multi-period inventory routing problem, including a supplier, set of customers, and a fleet of heterogeneous vehicles. Due to the complex nature of the IRP, we developed a Modified Adaptive Genetic Algorithm (MAGA) to solve a variety of instances efficiently. As a benchmark, we considered the results obtained by Cplex software and an efficient heuristic from the literature. Our approach showed significant improvement comparing the other two methods.
% \end{abstract}
% \begin{keywords}
% Inventory Routing Problem \sep Genetic Algorithm \sep Vendor Managed Inventory \sep Supply Chain Management \sep Adaptive Heuristic
% \end{keywords}

\let\WriteBookmarks\relax
\def\floatpagepagefraction{1}
\def\textpagefraction{.001}
\shorttitle{}
% \shortauthors{Meysam Mahjoob et~al.}

\title [mode = title]{A Modified Adaptive Genetic Algorithm for Multi-product Multi-period Inventory Routing Problem}                      
% \tnotemark[1,2]

% \tnotetext[1]{This document is the results of the research
%   project funded by the National Science Foundation.}

% \tnotetext[2]{The second title footnote which is a longer text matter
%   to fill through the whole text width and overflow into
%   another line in the footnotes area of the first page.}

\author[1]{Meysam Mahjoob}[
% type=editor,
%                         auid=000,bioid=1,
%                         prefix=Sir,
%                         role=Researcher,
%                         orcid=0000-0001-7511-2910
                        ]

\ead{mahjoob_m@alumni.ut.ac.ir}
% \ead[url]{www.cvr.cc, cvr@sayahna.org}

% \credit{Conceptualization of this study, Methodology, Software}

\address[1]{Department of Industrial Engineering, University of Tehran, Fooman, Rasht}
% \cormark[1]
\author[2]{Seyed Sajjad Fazeli}[]
\cormark[1]
\ead{sajjad.fazeli@wayne.edu}
% \ead[url]{www.cvr.cc, cvr@sayahna.org}

% \credit{Conceptualization of this study, Methodology, Software}

\address[2]{Department of Industrial and System Engineering, Wayne State university, Detroit, MI}

\author[3]{Soodabeh Milanlouei}[%
%   role=Co-ordinator,
%   suffix=Jr,
   ]
% \ead{cvr3@sayahna.org}
% \ead[URL]{www.sayahna.org}

% \credit{Data curation, Writing - Original draft preparation}

\ead{milanlouei.s@northeastern.edu}
% \ead[url]{www.cvr.cc, cvr@sayahna.org}

% \credit{Conceptualization of this study, Methodology, Software}

\address[3]{Center for Complex Network Research, Northeastern University, Boston, MA}

\author%
[4]
{Leyla Sadat Tavassoli}

\ead{Leylasadat.tavassoli@mavs.uta.edu}

\address[4]{Department of Industrial Manufacturing and Systems Engineering, University of Texas at Arlington, Arlington, TX}

\cortext[cor1]{Corresponding author: sajjad.fazeli@wayne.edu}
% \cortext[cor2]{Principal corresponding author}
% \fntext[fn1]{This is the first author footnote. but is common to third
% \fntext[fn2]{Another author footnote, this is a very long footnote and
%   it should be a really long footnote. But this footnote is not yet
%   sufficiently long enough to make two lines of footnote text.}

% \nonumnote{This note has no numbers. In this work we demonstrate $a_b$
%   the formation Y\_1 of a new type of polariton on the interface
%   between a cuprous oxide slab and a polystyrene micro-sphere placed

\begin{abstract}
Recent developments in urbanization and e-commerce have pushed businesses to deploy efficient systems to decrease their supply chain cost. Vendor Managed Inventory (VMI) is one of the most widely used strategies to effectively manage supply chains with multiple parties. VMI implementation asks for solving the Inventory Routing Problem (IRP). This study considers a multi-product multi-period inventory routing problem, including a supplier, set of customers, and a fleet of heterogeneous vehicles. Due to the complex nature of the IRP, we developed a Modified Adaptive Genetic Algorithm (MAGA) to solve a variety of instances efficiently. As a benchmark, we considered the results obtained by Cplex software and an efficient heuristic from the literature.  Through extensive computational experiments on a set of randomly generated instances, and using different metrics, we show that our approach distinctly outperforms the other two methods.

\end{abstract}

% \includegraphics{figs/grabs.pdf}

% \begin{highlights}
% \item Research highlights item 1
% \item Research highlights item 2
% \item Research highlights item 3
% \end{highlights}

\begin{keywords}
Inventory Routing Problem \sep Genetic Algorithm \sep Vendor Managed Inventory \sep Supply Chain Management \sep Adaptive Heuristic
\end{keywords}

\maketitle

\section{Introduction}
Throughout the past few decades, advancements in information and communications technology have led to novel and innovative businesses in supply chain management. Vendor Managed Inventory (VMI) is an example of such business models, representing an important paradigm in which the vendor (supplier) has full responsibility for controlling the retailer inventory using the inventory data that the retailer provides. Unlike traditional inventory management practices where retailers make their own decisions about the new orders' size and time, a vendor monitors the retailer’s inventory and acts as the decision-maker in the VMI model \cite{qin2014local}.
In this distribution model, deciding on the inventory level for both supplier and retailer depends on the timing and the number of deliveries that a retailer needs, which itself is influenced by the capacity of the vehicles used for delivery. Simultaneous decision making is a vital component to achieve cost-effectiveness in these systems. The main advantage of implementing a VMI system is that the vendor can increase the service level while reducing the distribution costs by using vehicles effectively. Retailers, on the other hand, can release the resources they usually use for inventory management. In VMI implementation, a vendor should solve an Integrated Routing Problem (IRP) to specify a distribution plan, minimizing the long-run average distribution and inventory costs throughout the supply chain, and avoiding shortage \cite{cordeau2015decomposition}.
The IRP simultaneously takes into account both vehicle routing and inventory management in a supply chain system and tries to maximize the overall performance. In the IRP, the supply chain system consists of a supplier and a number of retailers who are distributed geographically with certain demand levels. Within a specified time window, products are shipped from the supplier to retailers using a several vehicles. The demand of each customer is met by any vehicle per period.
In the VMI system, the supplier determines each retailer's replenishment policy and the vehicle routes to deliver the products while ensuring that no shortage occurs. The IRP attempts to minimize the overall logistics costs to meet retailers’ demands in a planning horizon. Under this framework, three main decisions need to be made: 1) when to fulfill the demand of a retailer; 2) the amount of products to be delivered to a retailer; 3) the optimal routes for vehicles. 
The IRP literature can be investigated through the modeling and methodology perspectives. From the modeling point of view, the literature can be classified into different groups based on the length of the planning horizon(single-period, multi-period, and infinite), type of product (single and multi), and delivery (non-split and split). One of the earliest studies was done by \cite{federgruen1984combined} where authors formulated the IRP as a single-period problem considering a plant and customers with uncertain demand. The research work by \cite{azadeh2017genetic} considered a single product IRP with transshipment where the products deteriorate during the time of storage or at warehouses. Authors in \cite{archetti2012hybrid} considered a multi-period single-product single-vehicle IRPwith a limited customers to minimize the inventory and transportation cost. An IRP with split delivery was studied by \cite{yu2008new}, where the constraint of the VRP that each customer should be served by exactly one vehicle was relaxed.
From the solution methodology perspective, the researchers have chosen different approaches toward solving the IRP problem. Although the IRP is NP-Hard, some researchers attempt to use exact methods to solve the IRP. For example, authors in \cite{coelho2013branch} developed a Branch-and-Cut (B\&C) algorithm to solve a multi-product multi-vehicle inventory-routing problem. A three-phase exact approach is developed by \cite{bertazzi2020exact} where in the first phase, a branch-and-cut algorithm is applied on extended IRP to obtain an extended lower-bound on cost of IRP. If the optimal integer solution could not be obtained in the first phase, the second phase is activated where a heuristic is applied to provide an upper-bound on the optimal cost. The upper bound will then be used in the third phase, where a branch-and-cut method is applied on the IRP. While most of the exact methods for IRP are based on the B\&C method, \cite{desaulniers2016branch} developed a branch-price algorithm by incorporating a group of valid inequalities, and a labeling algorithm
to solve column generation subproblems, combined with multiple acceleration techniques.
\par
Due to the complex nature of the IRP problem, many authors developed heuristic approaches to solve the IRP. A Variable Neighborhood Search (VNS) heuristic is developed by \cite{popovic2012variable} to solve a multi-period multi-product IRP. They implemented three different methods to generate the initial solutions. The authors in \cite{ramkumar2011hybrid} proposed a three-phase hybrid heuristic, including inventory allocation, warehouse clustering, and routing decisions. A Particle Swarm Optimization–Differential Evolution (PSO-DE) is presented by \cite{de2017sustainable} to solve green maritime IRP with time window constraints. They implemented a Genetic Algorithm (GA), Particle Swarm Optimization–Differential Evolution, and basic PSO to validate the results. A Genetic Algorithm is developed by \cite{moin2011efficient} to solve a many-to-one distribution network including a depot, an assembly plant and a number of suppliers. They suggested two types of modifications to their GA operations to tackle limitations on solving the problem. Authors in \cite{su2020matheuristic} proposed a metaheuristic algorithm to solve the special version of IRP. The algorithm follows a hierarchy-based consisting of a master problem for route adjustment and two subproblems for timing and flow decisions.
\par
In this study, we formulate a multi-product multi-period inventory routing problem (MMIRP) with a supplier, set of retailers (customers) distributed in a geographical area, and set of heterogeneous vehicles. We develop a Modified Adaptive Genetic Algorithm (MAGA) where we modify the crossover and mutation operations based on the chromosome representation matrix. Besides, to effectively define the GA parameters, we adopt an adaptive approach. The performance of the MAGA is tested against the Cplex software and a heuristic approach from the literature. We considered three different metrics to compare the results.
The remainder of this paper is organized as follows: Section \ref{formula} present problem definition, assumptions, and mathematical formulation. Section \ref{Method} introduces the solution methodology, where we provide details regarding the proposed heuristic. Section \ref{exper} presents instance characteristics, benchmark methods, and extensive computational experiments. Finally, Section \ref{Con} provides concluding remarks.

\section{Model Formulation} \label{formula}
\subsection{Problem Definition}
In this work, the MMIRP is defined as follows: a given set of vehicles which are dispatched from a supplier ($w_0$), a set of customers which are distributed in a geographical area; the goal is to determine the efficient routes for vehicles to satisfy the customers' demand while minimizing the fixed fleet cost, inventory, and transportation cost. The detailed assumptions regarding the model are provided as follows:
\begin{itemize}
    \item Different types of products are distributed in the network.
    \item There are no restrictions on the availability of products for the supplier. 
    \item Vehicles are heterogeneous, meaning that they can have different loading capacities.
    \item vehicles can visit multiple customers during their trip.
    \item The demand of customers is deterministic .
    \item Shortage is not allowed.
\end{itemize}
\subsection{Notation}
\begin{itemize}
	\item Sets
	\begin{itemize}
		\item $I$: Set of customers
		\item $\bar{I}$: Set of customers and supplier where $\bar{I}= I \cup \{w_0 \}$
		\item $T$: Set of time periods 
		\item $V$: Set of vehicles
	    \item $P$: Set of products types
	\end{itemize} 
	\item Model parameters
	\begin{itemize}
		\item $c_{i,j}$: Travel cost associated with edge $(i,j) \in \bar{I}$.
	\item $q^{v}$: Capacity of vehicle $v \in V$
	\item $f_{t}^{v}$: Fixed cost of vehicle $v \in V$ in period $t \in T$
	\item $s_i$: Storage capacity of customer $i \in I$
	\item $\alpha^p$: Weight associated with product $p \in P$
	\item $d^{p}_{i,t}$: Demand of customer $i \in I$ for product $p \in P$ in time $t \in T$
	\item $h^{p}_{i}$: Inventory cost of per unit product $p$ for customer $i$
	\end{itemize}
	\item  Decision variables
	\begin{itemize}
		\item $x_{i,j,t}^{v}$: 1 if the edge $(i, j) \in \bar{I}$ is traversed by vehicle $v \in V$ in period $t \in T$, and 0 otherwise;
		\item $y_{i,j,t}^{v,p}$: Amount of product $p \in P$ carried by vehicle $v \in V$ on edge $(i, j) \in \bar{I}$ .
        % \item $z^{p}_{i,t}$: Amount of delivery of product $p \in P$ to customer $i \in I$ in period $t \in T$
        \item $r_{i,t}^{p}$: Amount of on-hand inventory of product $p \in P$ by customer $i \in I$ 
	\end{itemize}
\end{itemize}
\subsection{Mathematical Formulation}\label{model}
\begin{alignat}{2}
& \hspace{-1.5cm}\text{Min} \sum_{t \in T} \bigg( \sum_{i \in I}\sum_{v \in V} f_{t}^{v} x^{v}_{w_{0},i,t} + \sum_{(i,j) \in \bar{I}}\sum_{v \in V} c_{i,j}x^{v}_{i,j,t} +\sum_{i \in I}\sum_{p \in P}h^{p}_{i}r_{i,t}^{p}\bigg)  && \label{obj} \\ 
&\hspace{-1.25cm} \text{s.t. } && \nonumber \\
&\hspace{-0.7cm} \sum_{j\in \bar{I}}x^{v}_{i,j,t}\leq 1 \qquad \forall   i \in \bar{I},t \in T, v \in V,  && \label{Cons2} \\ 
&\hspace{-0.7cm} \sum_{j \in \bar{I}}x^{v}_{i,j,t}- \sum_{k \in \bar{I}}x^{v}_{k,i,t}=0  \quad \forall   i \in \bar{I},t \in T, v \in V,  &&   \label{Cons3} \\ 
&\hspace{-0.7cm} \sum_{p \in P}\alpha^{p}y_{i,j,t}^{v,p}\leq q^{v}x_{i,j,t}^{v}\hspace{0.1cm} \forall (i,j) \in \bar{I}, t \in T, v \in V, && \label{Cons4} \\
& \hspace{-0.7cm}\sum_{j \in \bar{I}}\alpha^{p}y_{j,i,t}^{v,p}-\sum_{k \in \bar{I}}\alpha^{p}y_{i,k,t}^{v,p} \geq 0  \quad \forall  i \in \bar{I}, t \in T, v \in V, p\in P &&  \label{Cons5} \\ 
&  \hspace{-0.7cm} r_{i,t-1}^{p}-r_{i,t}^{p}+\sum_{v\in V}\bigg(\sum_{j \in \bar{I}}y_{j,i,t}^{v,p}-\sum_{k \in \bar{I}}y_{i,k,t}^{v,p}\bigg)=d^{p}_{i,t} &&\nonumber \\  
&\hspace{3.5cm} \quad \forall i\in \bar{I}, t \in T, p \in P \label{Cons6} \\
& \hspace{-0.7cm} \sum_{p \in P}\alpha^{p}r_{i,t}^{p}\leq s_{i} \qquad  \forall i \in I,t \in T,&&  \label{Cons7} \\
& \hspace{-0.7cm} x^{v}_{i,j,t}\in \{0,1\} \qquad \forall (i,j)\in \bar{I}, t \in T, v \in V,&&\label{Cons8}\\
& \hspace{-0.7cm} r_{i,t}^{p},y_{i,k,t}^{v,p}\geq 0 \qquad \forall i\in \bar{I}, p \in P, t \in T, v \in V&&\label{Cons9} 
\end{alignat}
The objective function (\ref{obj}) minimizes the fleet fixed, transportation, and inventory cost. Constraints (\ref{Cons2}) ensure that each customer is visited at most once in every period. Connectivity of routes for every vehicle is preserved by constraints \ref{Cons3}. Constraints (\ref{Cons4}) ensure that the amount of product carried by each vehicle does not violate the vehicle's capacity. The sub-tour in any feasible solution is eliminated by constraint (\ref{Cons5}). Constraints (\ref{Cons6}) guarantee the balance between the demand and on-hand inventory. Also, they ensure that the demand of customers is satisfied in each period The customers' storage capacity violation is prohibited by constraints (\ref{Cons7}). Finally, constraints (\ref{Cons8}) and (\ref{Cons9}) define the restrictions for the variables.
\section{Methodology and Algorithm Development}\label{Method}
The MMIRP could be solved by commercial solvers. However, to solve the mid and large-scale instances we need to develop algorithms to efficiently solve the MMIRP within a reasonable time limit. In the next section we comprehensively describe the MAGA that we developed to solve the MMIRP.
\subsection{Modified Genetic Algorithm}
The genetic algorithm is a stochastic optimization technique inspired by the process of natural selection, which is widely applied to solve different classes of NP-Hard problems \cite{ho2008hybrid,mahjoob2021green,tavana2018evolutionary,tavassoli2020integrated}. GA maintains a population of candidate solutions through the selective procedure. GA is initialized by a set of solutions called population. Each solution in the population is called a chromosome. Chromosomes progress through successive iterations called generations. Throughout each iteration, the chromosomes are evaluated by a fitness function. The fitter chromosomes have higher chances of being selected for GA operations such as mutation and crossover. The GA operations choose some parents and produce several offsprings as the new solutions. The new solutions are then accepted or rejected based on their fitness values as well as the solutions from the previous iterations to keep the population size fixed. The GA may converge to the best solution after a certain number of iterations. To solve the MMIRP using the GA algorithm, we have three major challenges: 1) The delivery schedule to meet customers' demand; 2) finding the best route for each vehicle; and 3) maintaining feasibility such that the vehicles and customer's capacity are not violated. Each of these is complex and difficult to solve, hence, a naive GA may not perform well for this problem. Therefore, to overcome the complexity, we propose a modified adaptive GA to efficiently solve the MMIRP.
\subsubsection{Solution Representation}
To represent the solution of the MMIRP, we define a chromosome in the form of a binary matrix of size ($|I| \times |T|$) where $|I|$ and $|T|$ denote the number of customers and periods, respectively. The elements $a_{i,t}=1$ indicates that the customer $i$ in period $t$ should be served. The representation matrix enables us to determine the amount of product for delivery to each customer in each period. To better demonstrate the advantage of the representation matrix, we provide an example in the next section.
\subsubsection{Explanatory Example}
\begin{figure*}[h]
% \centering
	\includegraphics[scale=0.45]{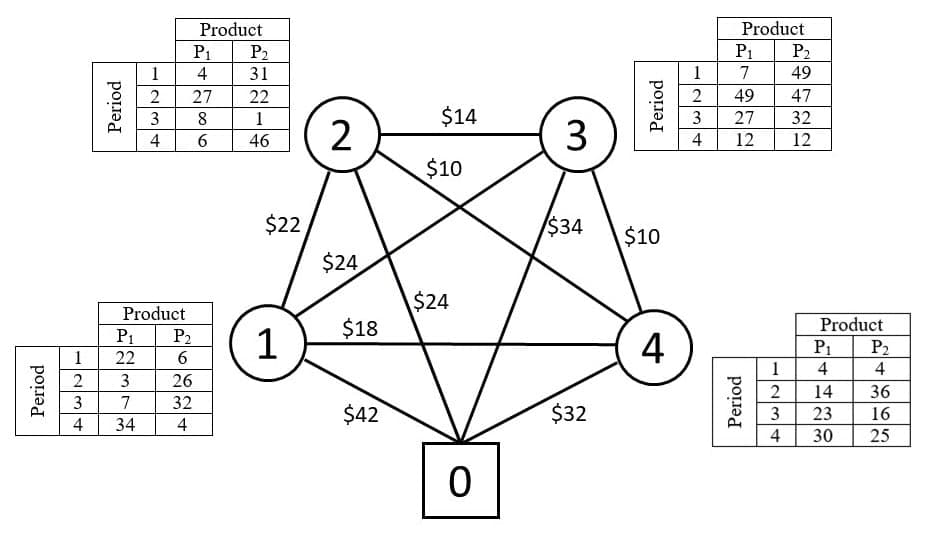}
% 	\vspace{-10mm}
    	\centering
    \captionsetup{justification=centering}

	\caption {Explanatory example network with a supplier and four customers}
	\label{Network}
\end{figure*}

\begin{figure*}
\centering
	\includegraphics[scale=0.5]{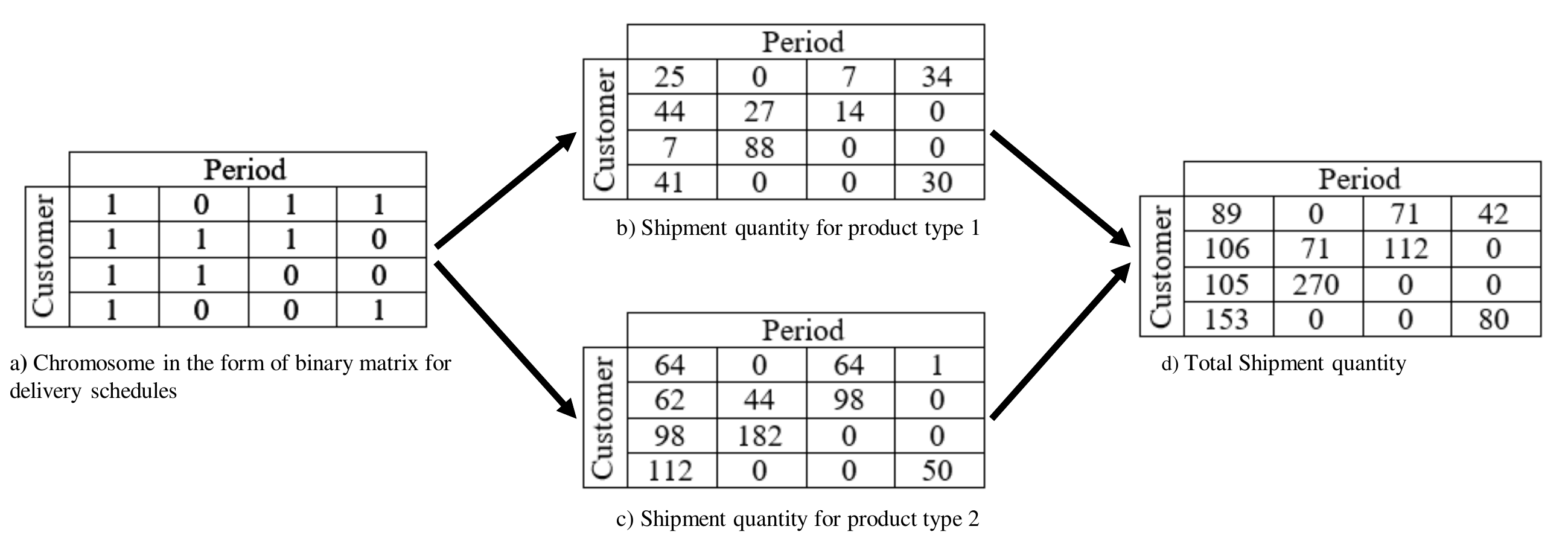}
% 	\vspace{-10mm}
    \captionsetup{justification=centering}
	\caption {A feasible solution for delivery schedules and shipment quantities represented in a matrix form}
	\label{feas-rep}
\end{figure*}

\begin{figure*}[h]
\centering
	\includegraphics[scale=0.5]{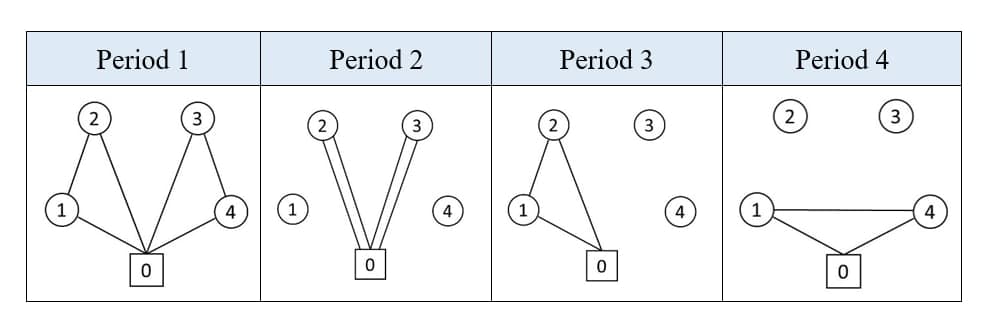}
% 	\vspace{-10mm}
    \captionsetup{justification=centering}
	\caption {Optimal or near-optimal routes for vehicles obtained by LKH method}
	\label{lkh}
\end{figure*}

We consider a network with a supplier, four customers, two types of product, four periods, and two vehicles with capacities of 300 and 400 units. Also, the fixed cost for each of the vehicles is equals to 10 units. The inventory costs of the two products for each customer in the distribution network are provided by mini-tables in Figure \ref{Network}. Also, the transportation cost for each route is shown on the edges. A feasible solution for this example is illustrated in Figure \ref{feas-rep}. Table (a) in this figure is the designed chromosome represented in the form of a binary matrix. Based on this matrix, customer 1 is served in periods 1, 3, and 4. Therefore, there is no delivery in period 2, and customer 2's demand is satisfied by the delivery in period 1. Similarly, the delivery schedule for the rest of the customers is shown in table (a) in Figure \ref{feas-rep}. Tables (b) and (c) in Figure \ref{feas-rep} indicate the shipment quantity based on the weight of products 1 and 2, respectively. For example, the shipment quantity of products 1 and 2 in period 1 is $1\times(22+3)$ and $2\times(6+26)$ respectively, where 1 and 2 are products' weight coefficients. Finally, table (d) \ref{feas-rep} shows the total shipment quantity. Once the shipment quantity for vehicles are calculated, we solve the single-vehicle routing problem by applying the Lin-Kernighan-Helsgaun heuristic (LKH) \cite{helsgaun2000effective}. Figure \ref{lkh} shows the optimal or near-optimal routes obtained by LKH.

As it is stated before, there are three types of cost in the distribution network: 1) fleet's fixed cost, 2) transportation cost, and 3) inventory cost. The first two costs are calculated based on the solution obtained in Figure \ref{lkh}. The inventory cost calculation depends on the period in which the products are delivered. For example, for the third customer, since the demand for periods 3 and 4 are delivered in period 2, we calculate the inventory cost as $1\times (27+39) + 2\times(12+12) = 114$ where the coefficients 1 and 2 indicate the number of periods that products are stored as inventory.\par
The feasibility of the chromosome is preserved considering the following conditions:
\begin{itemize}
    \item Since the shortage is not allowed, the representation matrix's elements in the first column should be 1.
    \item Split delivery should be avoided.
    \item In each period, the shipment product should not exceed vehicles' capacity.
    \item The amount of delivery to each customer minus the customer's demand at each period should not exceed the storage capacity. For instance, suppose that the storage capacity for customer 4 is violated in period 1, then we schedule a new delivery one period before the scheduled delivery ($a_{4,3}=1$) to avoid storage capacity violation.
\end{itemize}
So, every chromosome is checked based on the above conditions during the GA process, and the infeasible ones are corrected or dismissed from the solution pool.

\subsubsection{Selection}
 GA's convergence could be significantly affected by chromosome selection. The roulette wheel selection was first introduced by \cite{davis1985applying} to select the chromosomes for GA operations. Each section of the roulette wheel is assigned to a chromosome based on the magnitude of its fitness value. The fitness value of each chromosome is equal to the summation of fleet cost, transportation and inventory cost. The fitness values of the chromosomes determine their chance of being selected. We summarize the selection procedure in Algorithm \ref{Selection}. In the algorithm, the selection and cumulative probabilities of a chromosome $c$ are referred as $p_c$ and $g_c$, respectively, and the population size as $psize$. 
 \begin{algorithm}[!htbp]
	\caption{: Selection}
	\begin{algorithmic}
    	\State \text{\textit{Step 1:} Calculate the total fitness F:}
    	\State \hskip1.5em $$F = \sum_{c=1}^{psize} f_{c}$$
    	\State \text{\textit{Step 2:} Calculate the selection probability for each}
    	\State \hskip3em chromosome:
    	\State \hskip1.5em $$p_{c} = \frac{F-f_{c}}{F\times (psize-1)}$$
    	\State \text{\textit{Step 3:} Calculate the cumulative probability for each }
    	\State \hskip3em chromosome:
    	$$g_{c} =\sum_{i=1}^{c}p_{j} \quad c= 1,2,...,psize.$$
    	\State \text{\textit{Step 4:} Select the chromosome $c$ if :}
    	\State \hskip1.5em  $$g_{c-1}<r\leq g_{c} $$
    	\State \hskip1.5em where $r$ is a random number in range (0,1].
    \end{algorithmic}
    \label{Selection}
\end{algorithm}

\subsubsection{Crossover}
\begin{figure}
\centering
	\includegraphics[width=\columnwidth]{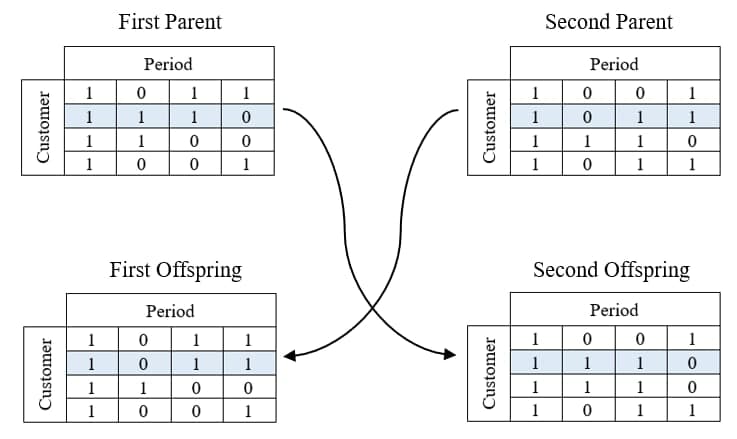}
% 	\vspace{-10mm}
	\caption {The crossover operation on matrix chromosome}
	\label{Crossover}
\end{figure}

The crossover operator is an essential part of the genetic algorithm, where two chromosomes (parents) contribute characteristics in creating a new chromosome (offspring), by randomly exchanging information. Crossover operator stochastically generates new solutions from an existing population. Crossover operation increases GA's ability in searching because it helps inherit and blend the good properties of the parents in the offspring from the elite solutions among the population. Since the designed chromosome in this work is constructed in a matrix form, the crossover operators can be applied on the rows or columns of the matrix. At each iteration, we randomly select a row or column from each parent and exchange the selected row or column among the parents. Figure (\ref{Crossover}) illustrates the crossover operator on the matrix chromosome.\\
\subsubsection{Mutation}
\begin{figure}
\centering
	\includegraphics[width=\columnwidth]{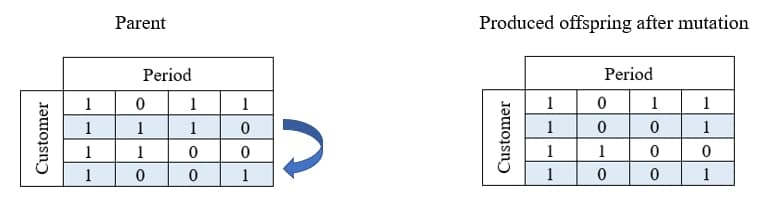}
% 	\vspace{-10mm}
	\caption {The mutation operation on matrix chromosome}
	\label{mutation}
\end{figure}

The mutation is another genetic operator that is used to explore new solutions in the solution space. In some newly formed offspring, some of their genes can be mutated with a low random probability and they may find characteristics that do not belong to any of their parents. If the mutation probability is too high, the searching process will transform into a primitive random search. The mutation occurs to maintain the diversity within the solution population and avoid trapping in local optima. The mutation operator designed for this problem is as follows: Once the roulette wheel algorithm selects a parent, we randomly select two rows or columns and flip their position. The new chromosome is then checked by the feasibility conditions and then added to the pool. Figure \ref{mutation} shows the mutation procedure.
\subsubsection{Adaptive GA Operators (AGO)}
The performance of GA highly depends on its parameters. In our study, we use the idea developed by \cite{mak2000adaptive}. In their approach, they used the parents' and offsprings' fitness values in each generation to create adaptive crossover and mutation. The mutation and crossover rates are adjusted based on the returned fitness value. The better the fitness value returned by any operator, the higher chance of operator rate. Therefore, the AGO strengthens the well-performing operators to generate more offsprings. We summarize the AGO steps in \ref{adap}. In this algorithm we denote parent size, offspring size, mutation rate, and crossover rate at iteration $k$ as $pr^{k}_{size}$, $of^{k}_{size}$, $mr^{k}$, $cr^{k}$. Also, $\bar{F}^{k}_{{pr}_{size}}$ and $\bar{F}^{k}_{{of}_{size}}$ indicate the average fitness value of parents and offsprings. We define a stopping criterion as $k_{max}$ which denotes the number of non-improving steps.
\begin{algorithm}[!htbp]
\caption{: AGO}
	\begin{algorithmic}
    	\State \textbf{Initialization:}
    	\State \hskip1.5em Set: \textit{$k$ $\gets$ $1$},
    	\State \textbf{while}  $k \leq k_{max}$ :
    	\State \hskip2em \textbf{if} $\frac{\bar{F}^{k}_{{pr}_{size}}}{\bar{F}^{k}_{{of}_{size}}}-1 \geq 0.1$:
    	\State \hskip3em  $mr^{k+1} \gets cr^{k}+0.005$, $cr^{k+1} \gets cr^{k}+0.05$
	    \State \hskip2em \textbf{elif} $\frac{\bar{F}^{k}_{{pr}_{size}}}{\bar{F}^{k}_{{of}_{size}}}-1 \leq 0.1$:
	    	\State \hskip3em  $mr^{k+1} \gets cr^{k}-0.005$, $cr^{k+1} \gets cr^{k}-0.05$
    	 \State \hskip2em \textbf{else}:
    	\State \hskip3em  $mr^{k+1} \gets cr^{k}$, $cr^{k+1} \gets cr^{k}$
    \end{algorithmic}
    \label{adap}
\end{algorithm}

%%%%%%%%%%%%%%%%%%%%%%%%%%%%%%%%%%%%%%%% EXPERIMENTS %%%%%%%%%%%%%%%%%%%%%%%%%%%%%
\vspace{-5mm}
\section{Experiments and Results} \label{exper}
In this section, we compare the computational performance of the proposed GA and Constructive Heuristic for Multi-product Inventory Routing Problem (CHMPIRP) proposed by \cite{dabiri2012constructive}. The authors employed three performance measures to evaluate the quality of their heuristic as follows.\\
\textbf{Difficulty}: is calculated as follows:\\
	   $$\text{Difficulty (\%)}= \frac{\text{Upper bound} - \text{Lower bound}}{\text{Upper bound}}$$
The lower difficulty percentage shows the tighter bound found by the solver within the stipulated time-limit. This measure is defined to better represent the other two measures.\\
\textbf{Closeness}: The closeness represents the distance between the heuristic's objective value (HBV) and  the best lower bound obtained by Cplex in percentage:\\
$$\text{Closeness (\%)}= \frac{\text{HBV} - \text{Lower bound}}{\text{HBV}}$$
Large values of Closeness indicate the weak lower bound or poor performance of heuristic.\\
\textbf{Saving}: The Saving indicates the distance between the best upper bound obtained by heuristic's objective value (HBV) and lower bound of Cplex in percentage:\\
    $$\text{Saving (\%)}= \frac{\text{Upper bound} - \text{HBV}}{\text{Upper bound}}$$

Large values of Saving indicates the higher quality of heuristic solution comparing to the Cplex solver.\\

They also provided the lower and upper bounds obtained by solving the model using Cplex 11.0 software package. In our study, all the experiments were implemented in MATLAB 2019a using a computer with an Intel \textregistered, CPU Core(TM) i7-9750H, 2.60 GHz, and 8GB RAM. In total 96 instances were generated in a square grid of size [20, 20], and the supplier at (10,10). The other details regarding instances are provided in Table \ref{charact}.
% \begin{table*}
%     \centering
% \includegraphics{charact.JPG}
%     \caption{Caption}
%     \label{charact}
% \end{table*}
\begin{table}[h]
\captionsetup{justification=centering}
\caption{Characteristics of the generated instances}
\centering
	\begin{tabular}{|c|c|c|} 
		\hline
		\multirow{2}{*}{Parameters} & \multicolumn{2}{c|}{Number of product types}
		\\ \cline{2-3}
			& $|P|=2$& $|P|=5$ \\
		\hline
 	    \multirow{1}{*}{$|I|$} &[5,10,20,30] & [5,10,20,30]\\
 	    \hline
 	     \multirow{1}{*}{$|T|$} &[5,7] & [5,7]\\
 	     \hline
 	     \multirow{1}{*}{$|V|$} &[3,5] & [3,5]\\
 	      \hline
 	     \multirow{1}{*}{$q^v$} &$|N|\times 50$ & $|N|\times 100$\\
 	      \hline
 	     \multirow{1}{*}{$s^i$} &300& 500\\
 	        \hline
 	     \multirow{1}{*}{$\alpha^p$} &[1,2]& [0.25,0.75,1,1.5,2.5]\\
 	     \hline
	\end{tabular}
	\label{charact}
\end{table}
\subsection{Results}
\begin{figure*}[h]
    \centering
    \includegraphics[scale=0.45]{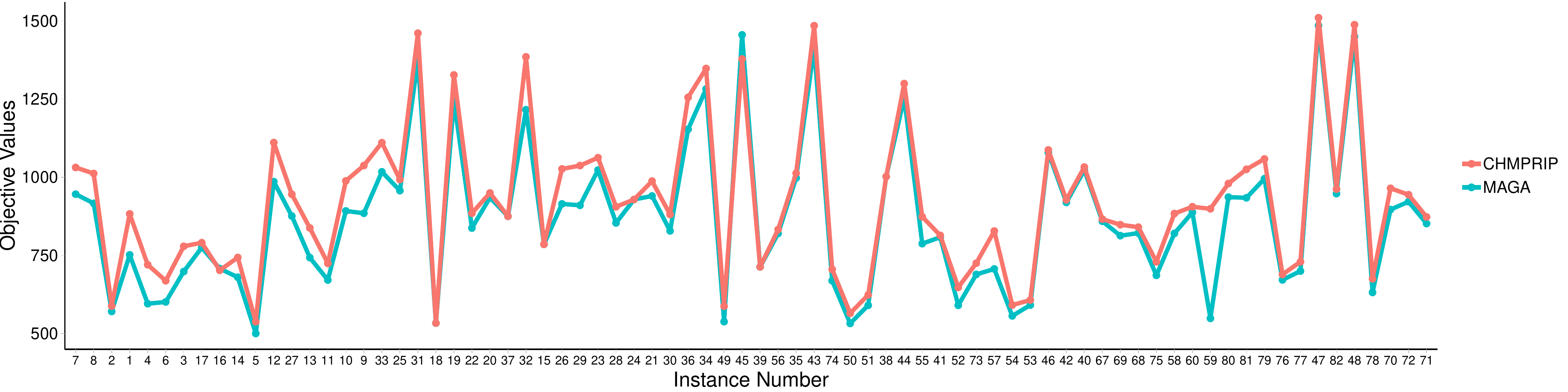}
    \captionsetup{justification=centering}
    \caption{Comparison between the MAGA and CHMPRIP in terms of the objective functions ranked by difficulty metric }
    \label{objs}
\end{figure*}

\begin{figure*}[h]
    \centering
    \captionsetup{justification=centering}
    \includegraphics[scale=0.45]{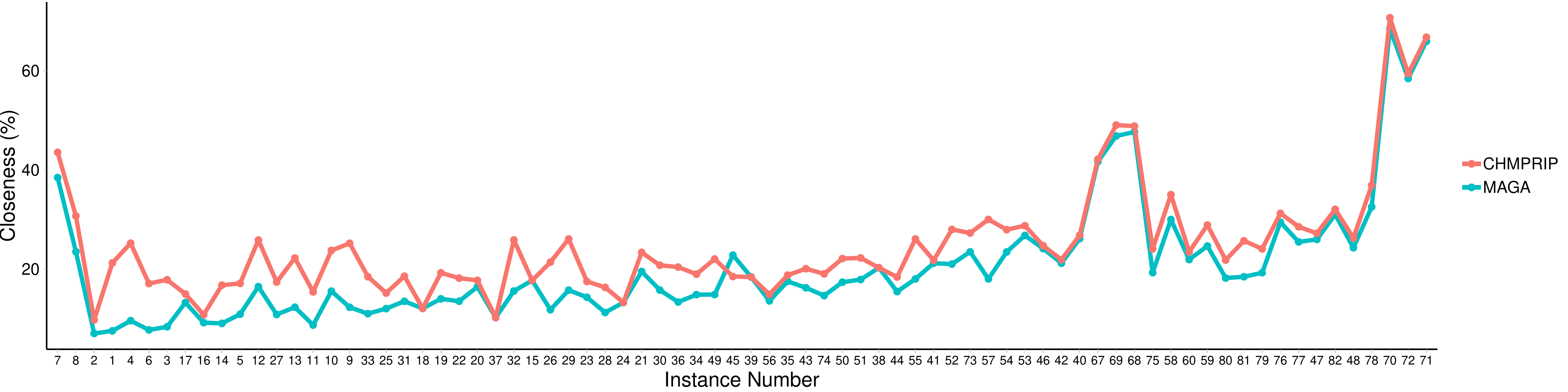}
    \caption{Comparison between the MAGA and CHMPRIP in terms of  closeness metric ranked by difficulty }
    \label{closeness}
\end{figure*}

\begin{figure*}[h]
    \centering
    \includegraphics[scale=0.45]{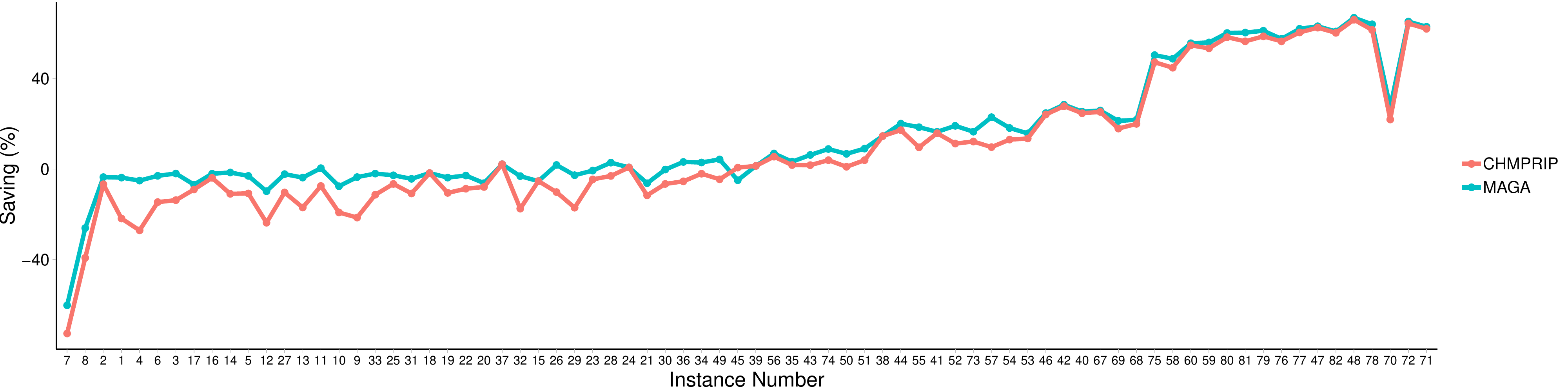}
    % \captionsetup{justification=raggedright,singlelinecheck=false}
    \captionsetup{justification=centering}
    \caption{Comparison between the MAGA and CHMPRIP in terms of  saving metric ranked by difficulty }
    \label{saving}
\end{figure*}

 Tables (\ref{results_1}) and (\ref{results_2}) compare the CHMPIRP with the MAGA. The first five columns indicate the instance number and the size of the instances in terms of the number of customers, periods, and vehicles. Columns 6 and 7 show the lower and upper bounds obtained by Cplex. The next two columns represent the best objective value returned by MAGA and CHMPIRP. Column 10 specifies the difficulty measure calculated by the lower and upper bounds obtained by Cplex. The closeness and saving measures calculated for MAGA and CHMPIRP are reported in the next four columns. The last two columns represent run time in seconds for MAGA and CHMPIRP, respectively. 
 \begin{table*}[h]
    \centering
    \captionsetup{justification=centering}
        \caption{Computational results for the generated instances}
        \vspace{-15mm}
\includegraphics[scale=0.9]{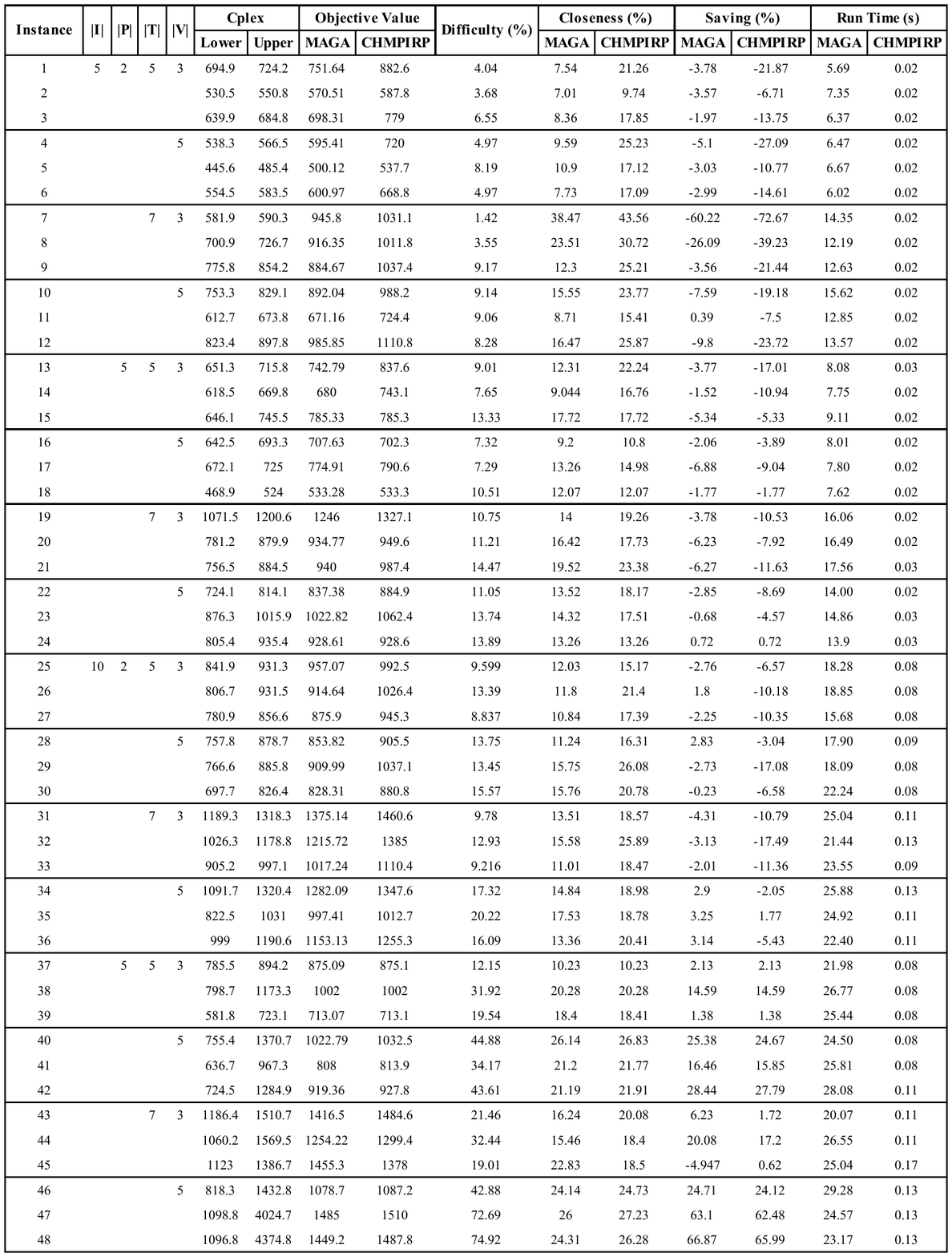}
    \label{results_1}
\end{table*}
 \begin{table*}[h]
    \centering
    \captionsetup{justification=centering}
    \caption{Computational results for the generated instances (continue)}
        \vspace{-15mm}
\includegraphics[scale=0.9]{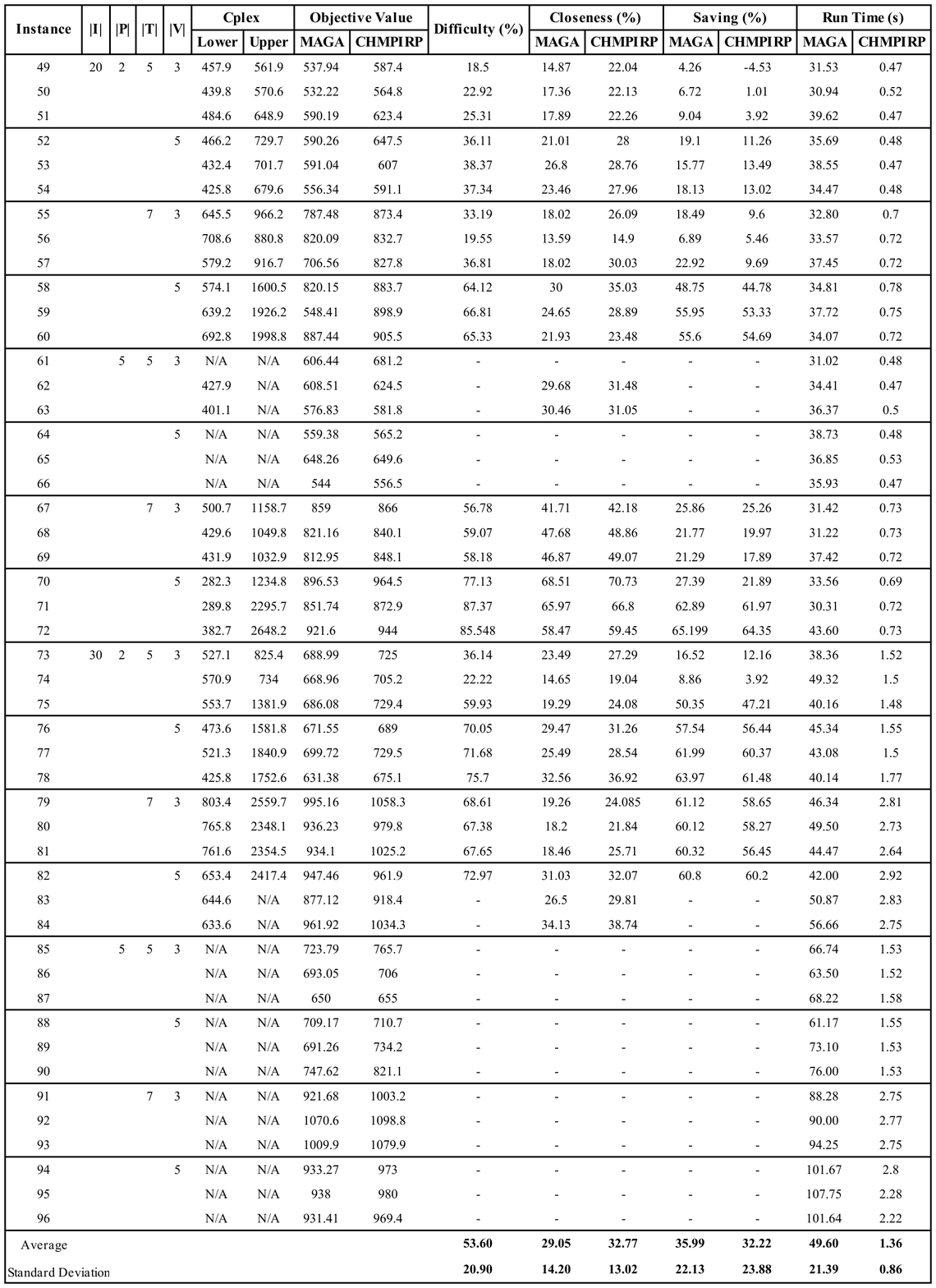}
    \label{results_2}
\end{table*}
%  \begin{table}
%     \centering
% \includegraphics{table3b.jpg}
%     \caption{Caption}
%     \label{results_2}
% \end{table}
In total, 96 instances were solved. In 16 instances the Cplex could not return any lower or upper bound within the stipulated time limit (2 hours). Also, in four instances, the Cplex could only return the lower bound. To better compare the MAGA's performance with CHMPRIP, we rank the problems based on the difficulty measure. Figure \ref{objs} illustrates the objective values returned by the two approaches ranked by difficulty percentage. Although the CHMPRIP offers a better performance in terms of run time, the results show that in almost all the instances, the MAGA outperformed the CHMPRIP by a significant margin in terms of the quality of solutions. Figures \ref{closeness} and \ref{saving} represent the performance of the two algorithms in terms of closeness and saving measures. The MAGA algorithm was able to provide better closeness in 95 instances by considerable difference comparing to the CHMPRIP. Also, in 94 instances, the MAGA indicated superior performance in terms of saving measure compared to the CHMPRIP.
% \newpage

% \clearpage
\subsection{Statistical Analysis}
To statistically compare the proposed method with the MAGA method, we use the paired t-test where the null hypothesis is that the mean difference between two sets of observations is zero. To do so, we solve each of the 96 problems five times and calculate the average for each problem separately for each method. We use a paired t-test to compare the means in these two populations, where observations in one method can be paired with observations in the other method. The obtained t-statistic is equal to -9.52, resulting in a p-value of \num{1.7e-15} which is smaller than the threshold of 0.05, indicating that there is strong evidence that the mean difference between the two methods is not equal to zero. With an average difference of -49.47 (95\% CI: -59.78, -39.15), the proposed method generally yields smaller objective values compared with the MAGA method.
\section{Conclusion}\label{Con}
In this research, we considered a multi-product multi-period inventory routing problem with a heterogeneous fleet of vehicles. We proposed a modified adaptive genetic algorithm to efficiently solve a variety of instances. The approach incorporated different techniques within the genetic algorithm framework to enhance the performance of the algorithm. Numerical studies were performed on randomly generated data sets from literature. The efficacy of the proposed method is benchmarked by comparing it to the commercial solver Cplex as well as a heuristic from the literature (CHMPRIP). The quality of the solutions is justified based on three different metrics. Other than a few instances, the results indicate that MAGA has outperformed the CHMPRIP in all metrics by a significant margin. Besides, we performed a statistical analysis to show the significant difference between the two heuristics in terms of solution quality. Future work could include incorporating green approaches by either replacing the vehicles with Electric vehicles or incorporating $\text{CO}_2$ emissions into the objective functions. Another extension could be incorporating uncertainties into the customers' demand. From an algorithmic perspective, the use of decomposition algorithms like the column generation techniques can also be evaluated.  

% \section{Ref}
\bibliographystyle{unsrt}
\bibliography{Refs}

\begin{thebibliography}{10}

\bibitem{qin2014local}
Lei Qin, Lixin Miao, Qingfang Ruan, and Ying Zhang.
\newblock A local search method for periodic inventory routing problem.
\newblock {\em Expert Systems with Applications}, 41(2):765--778, 2014.

\bibitem{cordeau2015decomposition}
Jean-Fran{\c{c}}ois Cordeau, Demetrio Lagan{\`a}, Roberto Musmanno, and
  Francesca Vocaturo.
\newblock A decomposition-based heuristic for the multiple-product
  inventory-routing problem.
\newblock {\em Computers \& Operations Research}, 55:153--166, 2015.

\bibitem{federgruen1984combined}
Awi Federgruen and Paul Zipkin.
\newblock A combined vehicle routing and inventory allocation problem.
\newblock {\em Operations research}, 32(5):1019--1037, 1984.

\bibitem{azadeh2017genetic}
A~Azadeh, S~Elahi, M~Hosseinabadi Farahani, and B~Nasirian.
\newblock A genetic algorithm-taguchi based approach to inventory routing
  problem of a single perishable product with transshipment.
\newblock {\em Computers \& Industrial Engineering}, 104:124--133, 2017.

\bibitem{archetti2012hybrid}
Claudia Archetti, Luca Bertazzi, Alain Hertz, and M~Grazia Speranza.
\newblock A hybrid heuristic for an inventory routing problem.
\newblock {\em INFORMS Journal on Computing}, 24(1):101--116, 2012.

\bibitem{yu2008new}
Yugang Yu, Haoxun Chen, and Feng Chu.
\newblock A new model and hybrid approach for large scale inventory routing
  problems.
\newblock {\em European Journal of Operational Research}, 189(3):1022--1040,
  2008.

\bibitem{coelho2013branch}
Leandro~C Coelho and Gilbert Laporte.
\newblock A branch-and-cut algorithm for the multi-product multi-vehicle
  inventory-routing problem.
\newblock {\em International Journal of Production Research},
  51(23-24):7156--7169, 2013.

\bibitem{bertazzi2020exact}
Luca Bertazzi, Demetrio Lagan{\`a}, Jeffrey~W Ohlmann, and Rosario Paradiso.
\newblock An exact approach for cyclic inbound inventory routing in a level
  production system.
\newblock {\em European Journal of Operational Research}, 283(3):915--928,
  2020.

\bibitem{desaulniers2016branch}
Guy Desaulniers, J{\o}rgen~G Rakke, and Leandro~C Coelho.
\newblock A branch-price-and-cut algorithm for the inventory-routing problem.
\newblock {\em Transportation Science}, 50(3):1060--1076, 2016.

\bibitem{popovic2012variable}
Dra{\v{z}}en Popovi{\'c}, Milorad Vidovi{\'c}, and Gordana Radivojevi{\'c}.
\newblock Variable neighborhood search heuristic for the inventory routing
  problem in fuel delivery.
\newblock {\em Expert Systems with Applications}, 39(18):13390--13398, 2012.

\bibitem{ramkumar2011hybrid}
N~Ramkumar, P~Subramanian, TT~Narendran, and K~Ganesh.
\newblock A hybrid heuristic for inventory routing problem.
\newblock {\em International Journal of Electronic Transport}, 1(1):45--63,
  2011.

\bibitem{de2017sustainable}
Arijit De, Sri~Krishna Kumar, Angappa Gunasekaran, and Manoj~Kumar Tiwari.
\newblock Sustainable maritime inventory routing problem with time window
  constraints.
\newblock {\em Engineering Applications of Artificial Intelligence}, 61:77--95,
  2017.

\bibitem{moin2011efficient}
Noor~Hasnah Moin, Said Salhi, and NAB Aziz.
\newblock An efficient hybrid genetic algorithm for the multi-product
  multi-period inventory routing problem.
\newblock {\em International Journal of Production Economics}, 133(1):334--343,
  2011.

\bibitem{su2020matheuristic}
Zhouxing Su, Zhipeng L{\"u}, Zhuo Wang, Yanmin Qi, and Una Benlic.
\newblock A matheuristic algorithm for the inventory routing problem.
\newblock {\em Transportation Science}, 54(2):330--354, 2020.

\bibitem{ho2008hybrid}
William Ho, George~TS Ho, Ping Ji, and Henry~CW Lau.
\newblock A hybrid genetic algorithm for the multi-depot vehicle routing
  problem.
\newblock {\em Eng. Appl. Artif. Intell.}, 21(4):548--557, 2008.

\bibitem{mahjoob2021green}
Meysam Mahjoob, Seyed~Sajjad Fazeli, Soodabeh Milanlouei, Ali~Kamali
  Mohammadzadeh, and Leyla~Sadat Tavassoli.
\newblock Green supply chain network design with emphasis on inventory
  decisions.
\newblock {\em arXiv preprint arXiv:2104.05924}, 2021.

\bibitem{tavana2018evolutionary}
Madjid Tavana, Kaveh Khalili-Damghani, Debora Di~Caprio, and Zeynab Oveisi.
\newblock An evolutionary computation approach to solving repairable
  multi-state multi-objective redundancy allocation problems.
\newblock {\em Neural Computing and Applications}, 30(1):127--139, 2018.

\bibitem{tavassoli2020integrated}
Leyla~Sadat Tavassoli, Nahal Sakhavand, and Seyed~Sajjad Fazeli.
\newblock Integrated preventive maintenance scheduling model with redundancy
  for cutting tools on a single machine.
\newblock {\em Engineering, Technology \& Applied Science Research},
  10(6):6542--6548, 2020.

\bibitem{helsgaun2000effective}
Keld Helsgaun.
\newblock An effective implementation of the lin--kernighan traveling salesman
  heuristic.
\newblock {\em Eur. J. Oper. Res.}, 126(1):106--130, 2000.

\bibitem{davis1985applying}
Lawrence Davis.
\newblock Applying adaptive algorithms to epistatic domains.
\newblock In {\em IJCAI}, volume~85, pages 162--164, 1985.

\bibitem{mak2000adaptive}
KL~Mak, YS~Wong, and XX~Wang.
\newblock An adaptive genetic algorithm for manufacturing cell formation.
\newblock {\em The International Journal of Advanced Manufacturing Technology},
  16(7):491--497, 2000.

\bibitem{dabiri2012constructive}
Nooraddin Dabiri, Mohammad~Jafar Tarokh, and Mostafa Setak.
\newblock A constructive heuristicfor amulti-productinventoryrouting problem.
\newblock {\em Transportation Research}, 2(1):11, 2012.

\end{thebibliography}

\end{document}